\documentclass[a4paper,12pt]{amsart}

\usepackage{array}

\usepackage{epsf}

\newcommand{\inspic}[1]{\begin{tabular}{c}\epsfbox{#1}\end{tabular}}

\newtheorem{theorem}{Theorem}
\newtheorem{lemma}{Lemma}

\def\oM{{\overline{\mathcal{M}}}}

\def\C{\mathbb{C}}
\def\Z{\mathbb{Z}}

\def\d{\partial}

\title{Belorousski-Pandharipande relation in dGBV algebras}

\thanks{S.S. was partially supported by grants RFBR-04-02-17227-a, NSh-1972.2003.1, NWO-RFBR-047.011.2004.026 (RFBR-05-02-89000-NWO-a), by the G\"oran Gus\-tafs\-son foundation, and by the P.~Deligne's fund based on his 2004 Balzan prize. I.S. was partially supported by grant RFBR-03-01-00167-a.}



\author{Sergei Shadrin}

\address{Department of Mathematics, Stockholm University, Stockholm, SE-10691, Sweden and Moscow Center for Continuous Mathematical Education, Boljshoi Vlasjevskii Pereulok 11, Moscow, 119002, Russia.}

\email{shadrin@math.su.se and shadrin@mccme.ru}

\author{Igor Shneiberg}

\address{Department of Algebra, Faculty of Mechanics and Mathematics, Moscow State University,
Leninskie Gory, GSP, Moscow, 119899, Russia.}

\email{shneiberg@mtu-net.ru}

\begin{document}

\begin{abstract}
We prove that the genus expansion of solutions of the WDVV equation constructed from dGBV algebras satisfy the differential equation determined by the Belorousski-Pandharipande relation in cohomology of the moduli space of curves $\oM_{2,3}$.
\end{abstract}

\maketitle


\section{Introduction}

In papers~\cite{ls,s}, the algebraic formalization of Zwiebach invariants gives a purely algebraic construction of genus expansion of solutions of the WDVV equation in terms of cH-algebras. Since we know that Zwiebach invarinats induce Gromov-Witten invariants on subbicomlexes with zero differential, then it is a natural question to check different relation coming from geometry of the moduli space of curves.

In Gromov-Witten theory, we represent a solution of the WDVV equation as a generating function for genus zero Gromov-Witten invariants of a siutable algebraic variety (without $\psi$-classes). Then we consider Gromov-Witten invariants with $\psi$-classes and in arbitrary genus. In this approach, any relation among natural strata in cohomology of the moduli space of curves gives us a differential equation for the Gromov-Witten potential. This property is just a corollary of the splitting axiom.

In our construction, we do not have a definition that can be compared with full Gromov-Witten potential. We just define genus expansion with descendants ($\psi$-classes) only at one point.
But some relations from Gromov-Witten theory are already nontrivial even for this reduced genus expansion. In particular, it is enough to have descendants at one point to pose a question on Belorousski-Pandharipande relation among codimension $2$ strata in $\oM_{2,3}$~\cite{bp}.

In this paper, we prove that our genus expansion satisfy the differential equation defined by the Belorousski-Pandharipande relation. In fact, it is the unique relation in genera $\leq 2$ that makes sense in our construction and that we have not yet checked in~\cite{ls,s}.

Since this paper is just a sequel to~\cite{ls,s}, we refer to these papers for the origin, motivation, and more detailed exposition of the new construction that we study here. Also we mention that in this paper we keep all ealier mathematical problems of this theory such as the lack of examples and the problem of convergence.

This paper is organized as follows. In Section~2, we recall our construction of genus expansion (in fact, only its parts used in this paper). In Section~3, we recall the Belorousski-Pandharipande relation and formulate our main theorem. In the rest of the paper, we prove (or rather outline the proof of) our theorem.

We are very grateful to A.~S.~Losev for the numerous of helpful remarks and discussions.

\section{Construction of potential}

In our case, the Belorousski-Pandharipande relation is a differential equation for four different formal power series: generating functions for the correlators in genera $0$, $1$, and $2$ without descendants ($\Phi_0$, $\Phi_1$, and $\Phi_2$, respectively), and generating function for the correlators in genus $2$ with one descendant at one point ($\Phi_2^{(1)}$). Our goal in this section is to define these four formal power series.

In subsection~\ref{cH-algebras} we explain what is a cH-algebra and fix notations for all necessary operators in it. In subsection~\ref{tensors} we explain how we use graph to encode tensor expressions. In subsection~\ref{usage} we fix notation for all tensors in cH-algebras that we use in this paper. Also we discuss there a subtlety related to the signs. Then in subsection~\ref{def-graphs} we define $\Phi_0$, $\Phi_1$, $\Phi_2$, and $\Phi_2^{(1)}$, and give precise formulas (in terms of graphs) for the first several terms of these power series.

In fact, it can also be useful to study the definition of the full potential with descendants only at one point given in~\cite{s}. But the full definition given there is rather involved and appears to be natural only in the course of calculations in~\cite{s} or in the framework of Zwiebach invariants~\cite{ls}.

\subsection{cH-algebras}\label{cH-algebras}
In this section, we recall the definition of cH-alge\-bras~\cite{ls,s}. A supercommutative associative $\C$-algebra $H$ is called cH-algebra, if there are two odd linear operators
$Q,G_-\colon H\to H$ and an integral $\int\colon H\to\C$ satisfying the following axioms:
\begin{enumerate}
\item $Q^2=G_-^2=QG_-+G_-Q=0$;
\item $H=H_0\oplus H_4$, where $QH_0=G_-H_0=0$ and $H_4$ is represented as a direct sum of subspaces of dimension $4$ generated by $e_\alpha, Qe_\alpha, G_-e_\alpha, QG_-e_\alpha$ for some vectors $e\in H_4$, i.~e. $H_4=\bigoplus_{\alpha}\, \langle e_\alpha, Qe_\alpha, G_-e_\alpha, QG_-e_\alpha \rangle$ (Hodge decomposition);
\item $Q$ is an operator of the first order, it satisfies the Leibniz rule: $Q(ab)=Q(a)b+(-1)^{\tilde a}aQ(b)$
(here and below we denote by $\tilde a$  the parity of $a\in H$);
\item $G_-$ is an operator of the second order, it satisfies the $7$-term relation:
$G_-(abc)=G_-(ab)c+(-1)^{\tilde b(\tilde a+1)}bG_-(ac)
+(-1)^{\tilde a}aG_-(bc)\\
-G_-(a)bc-(-1)^{\tilde a}aG_-(b)c
-(-1)^{\tilde a+\tilde b}abG_-(c).$
\item $G_-$ satisfies the property called $1/12$-axiom: $str(G_-\circ a\cdot)=(1/12)str(G_-(a)\cdot)$ (here $a\cdot$ and $G_-(a)\cdot$ are the operators of multiplication by $a$ and $G_-(a)$ respectively).
\end{enumerate}

We define an operator $G_+\colon H\to H$. We put $G_+H_0=0$; on each subspace $\langle e_\alpha, Qe_\alpha, G_-e_\alpha, QG_-e_\alpha \rangle$, we define $G_+$ as $G_+e_\alpha=G_+G_-e_\alpha=0$, $G_+Qe_\alpha=e_\alpha$, and $G_+QG_-e_\alpha =G_-e_\alpha$. We see that $[G_-,G_+]=0$;
$\Pi_4=[Q,G_+]$ is the projection to $H_4$ along $H_0$; $\Pi_0=\mathrm{Id}-\Pi_4$ is the projection to $H_0$ along $H_4$.

An integral on $H$ is an even linear function $\int\colon H\to\C$. We require
$\int Q(a)b =  (-1)^{\tilde a+1}\int aQ(b)$, $\int G_-(a)b = (-1)^{\tilde a}\int aG_-(b)$, and
$\int G_+(a)b = (-1)^{\tilde a}\int aG_+(b)$. These properties imply that $\int G_-G_+(a)b=\int aG_-G_+(b)$, $\int \Pi_4(a)b=\int a\Pi_4(b)$, and $\int \Pi_0(a)b=\int a\Pi_0(b)$.

We can define a scalar product on $H$: $(a,b)=\int ab.$ We suppose that this scalar product is non-degenerate. Using  scalar product we may turn an operator $A: H \to H$ into bivector that we
denote by $[A]$.

\subsection{Tensor expressions in terms of graphs}\label{tensors}
Here we explain a way to encode some tensor expressions over an arbitrary vector space in terms of graphs.

Consider an arbitrary graph (we allow graphs to have leaves and we require vertices to be at least of degree $3$). We associate a symmetric $n$-form to each internal vertex of degree $n$, a symmetric bivector to each egde, and a vector to each leaf. Then we can substitute the tensor product of all vectors in leaves and bivectors in edges into the product of $n$-forms in vertices, distributing the components of tensors in the same way as the corresponding edges and leaves are attached to vertices in the graph. This way we get a number.

Let us study an example:
\begin{equation}
\inspic{arb.1}
\end{equation}
We assign a $5$-form $x$ to the left vertex of this graph and a $3$-form $y$ to the right vertex. Then the number that we get from this graph is $x(a,b,c,v,w)\cdot y(v,w,d)$.

Note that vectors, bivectors and $n$-forms used in this construction can depend on some variables. Then what we get is not a number, but a function.

\subsection{Usage of graphs in cH-algebras}\label{usage}

Consider a cH-algebra $H$. There are some standard tensors over $H$, which we associate to elements of graphs below. Here we introduce the notations for these tensors.

We always assign the form
\begin{equation}
(a_1,\dots,a_n)\mapsto \int a_1\cdot\dots\cdot a_n
\end{equation}
to a vertex of degree $n$.

There is a collection of bivectors that will arise below at edges: $[G_-G_+]$, $[\Pi_0]$, $[Id]$, $[QG_+]$, $[G_+Q]$, $[G_+]$, and $[G_-]$. In pictures, edges with these bivectors will be denoted by
\begin{equation}
\inspic{arb.2},\quad \inspic{arb.3},\quad \inspic{arb.4},\quad \inspic{arb.5},\quad \inspic{arb.6},\quad \inspic{arb.7},\quad \inspic{arb.8},
\end{equation}
respectively. Note that an empty edge corresponding to the bivector $[Id]$ can usually be contracted (if it is not a loop).

The vectors that we will put at leaves depend on some variables. Let $\{e_1,\dots,e_s\}$ be a
homogeneous basis of $H_0$. To each vector $e_i$ we associate two formal variables, $T_{0,i}$ and $T_{1,i}$, of the same parity as $e_i$. Then we will put at a leaf either the vector $E_0=\sum e_i T_{0,i}$ (denoted by an empty leaf) or the vector $E_1=\sum e_i T_{1,i}$ (denoted by an arrow at the leaf).

\subsubsection{Remark}

There is a subtlety related to the fact that $H$ is a $\Z_2$-graded space. In order to give an honest definition we must do the following. Suppose we consider a graph of genus $g$. We can choose $g$ edges in such a way that the graph being cut at these edge turns into a tree. To each of these edges we have already assigned a bivector $[A]$ for some operator $A\colon H\to H$. Now we have to put the bivector $[JA]$ instead of the bivector $[A]$, where $J$ is an operator defined by the formula $J\colon a\mapsto (-1)^{\tilde a} a$.

In particular, consider the following graph (this is also an example to the notations given above):
\begin{equation}
\inspic{arb.9}
\end{equation}
An empty loop corresponds to the bivector $[Id]$. An empty leaf corresponds to the vector $E_0$. A trivalent vertex corresponds to the $3$-form given by the formula $(a,b,c)\mapsto\int abc$.

If we ignore this remark, then what we get is just the trace of the operator $a\mapsto E_0\cdot a$. But using this remark we get the supertrace of this operator.

In fact, this subtlety will play no role in this paper. It affects only some signs in calculations and all these signs will be hidden in lemmas shared from~\cite{ls,s}. So, one can just ignore this remark.

\subsection{Construction of $\Phi_0$, $\Phi_1$, $\Phi_2$, and $\Phi_2^{(1)}$}\label{def-graphs}

Now we describe $\Phi_0$, $\Phi_1$, $\Phi_2$, and $\Phi_2^{(1)}$ using the notations given above.

The formal power series $\Phi_0$ ($\Phi_1$, $\Phi_2$) is just the sum over all trivalent graphs of genus $0$ ($1$, $2$, respectively) with empty leaves and edges with thick black dots. At each graph we put a coefficient equals to the inverse of the number of its automorphisms.
\begin{align}
\Phi_0 & =\frac{1}{6}\inspic{potentia.15}+\frac{1}{8}\inspic{potentia.16} +\frac{1}{8}\inspic{potentia.17} +\dots\\
\Phi_1 & =\frac{1}{2}\inspic{potentia.18} +\frac{1}{4}\inspic{potentia.19} +\frac{1}{4} \inspic{potentia.20} +\dots\\
\Phi_2 &=\frac{1}{12}\inspic{potentia.1}+\frac{1}{8}\inspic{potentia.2}
+\frac{1}{8}\inspic{potentia.3}+\frac{1}{4}\inspic{potentia.4}
+\frac{1}{4}\inspic{potentia.5}+\dots
\end{align}

The formal power series $\Phi_2^{(1)}$ is the some over graphs of genus $2$ with edges with thick black dots satisfying some additional conditions. First, there is exactly one vertex of degree $4$ and all othe vertices are trivalent. Second, at this vertex of degree $4$, there is a leaf with arrow. Third, all other leaves are empty. Each graph is weighted with the inverse of the number of its automorphisms fixing the leaf with arrow.
\begin{multline}
\Phi_2^{(1)}= \frac{1}{12}\inspic{potentia.6}+\frac{1}{4}\inspic{potentia.7}
+\frac{1}{8}\inspic{potentia.8}+\frac{1}{4}\inspic{potentia.9}
+\frac{1}{4}\inspic{potentia.10} \\
+\frac{1}{4}\inspic{potentia.11}+\frac{1}{2}\inspic{potentia.12}+
\frac{1}{4}\inspic{potentia.13}+\frac{1}{4}\inspic{potentia.14}+\dots
\end{multline}
In fact, in order to obtain an expression for $\Phi_2^{(1)}$ one can just take the expression for $\Phi_2$ and add an additional leaf with arrow in all possible ways to each graph.
Also it is obvious that $\Phi_2^{(1)}$ is linear in $T_{1,i}$, $i=1,\dots,s$.


\section{Belorousski-Pandharipande relation}

\subsection{Notations}

The Belorousski-Pandharipande relation is a relation in (co)homology of $\oM_{2,3}$ between the cycles of natural strata of complex codimension $2$ in $\oM_{2,3}$. Below,
we list the strata participating in the Belorousski-Pandharipande relation:
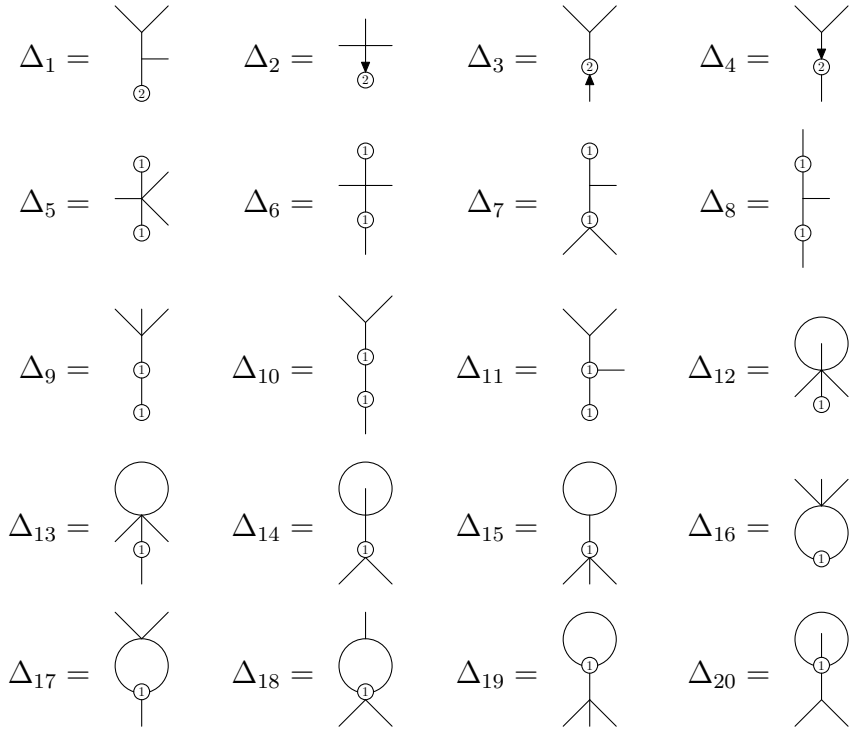
\begin{align*}
\Delta_{1} &=\inspic{bp.1} & \Delta_{2} &=\inspic{bp.2} & \Delta_{3} &=\inspic{bp.3} & \Delta_{4} &=\inspic{bp.4} \\
\Delta_{5} &=\inspic{bp.5} & \Delta_{6} &=\inspic{bp.6} & \Delta_{7} &=\inspic{bp.7} & \Delta_{8} &=\inspic{bp.8} \\
\Delta_{9} &=\inspic{bp.9} & \Delta_{10} &=\inspic{bp.10} &\Delta_{11} &=\inspic{bp.11} & \Delta_{12} &=\inspic{bp.12} \\
\Delta_{13} &=\inspic{bp.13} & \Delta_{14} &=\inspic{bp.14} &\Delta_{15} &=\inspic{bp.15} & \Delta_{16} &=\inspic{bp.16} \\
\Delta_{17} &=\inspic{bp.17} & \Delta_{18} &=\inspic{bp.18} &\Delta_{19} &=\inspic{bp.19} & \Delta_{20} &=\inspic{bp.20}
\end{align*}

We explain our notations. Note that graphs here have completely different meaning then all other graphs in this paper. We use the language of dual graphs, that is,
vertices correspond to irreducible curves,
edges corespond to points of intersection,
leaves correspond to marked points. A thick vertex
labeled by $2$ corresponds to a genus two curve;
a thick vertex labeled by $1$ corresponds to a genus one curve.
A simple vertex corresponds to a genus zero curve.
An arrow on an edge or a leaf means
that we take the $\psi$-class at the destination of the arrow.

This way to describe strata in the moduli space of curves was introduced by E.~Getzler, see~\cite{g}.

For example, consider the picture of $\Delta_1$. A generic point of this stratum is represented by a three component curve such that one component has genus $0$; there is one marked point on these curve and two other curves intersect this one. One of the other curves has genus $0$ and two marked point; another curve has genus $2$ and no marked points.

Another example. We consider the picture of $\Delta_2$. A generic point of this stratum is represented by a two-component curve; one component has genus zero; there are three marked point and one point of intersection with another curve. Another curve has genus $2$; there are no marked points, but we take the $\psi$-class on this curve at the point of intersection.

One more example is the picture of $\Delta_3$. A generic point of this stratum is represented by a two-component curve; one curve has genus $0$, two marked points, and one point of intersection with another curve. Another curve has genus $2$, one point of intersection with the first curve, and one marked point with $\psi$-class.

\subsection{The relation} We recall the Belorousski-Pandharipande relation:
\begin{align}\label{BP-relation}
-4\Delta_{1}+12\Delta_{2}+6\Delta_{3}-6\Delta_{4}+\frac{12}{5}\Delta_{5}
-\frac{12}{5}\Delta_{6}+\frac{24}{5}\Delta_{7} &\\
-\frac{36}{5}\Delta_{8}
-\frac{36}{5}\Delta_{9}+
\frac{18}{5}\Delta_{10}
-\frac{12}{5}\Delta_{11}
+\frac{1}{10}\Delta_{12}-
\frac{3}{10}\Delta_{13}&\notag\\
+\frac{3}{10}\Delta_{14}
-\frac{1}{10}\Delta_{15}
+\frac{6}{5}\Delta_{16}-\frac{6}{5}\Delta_{17}+\frac{2}{5}\Delta_{18}-
\frac{3}{5}\Delta_{19}-\frac{1}{5}\Delta_{20}& =0.\notag
\end{align}

One can notice that the coefficients in Equation~(\ref{BP-relation})
do not coinside with the coefficients of the initial relation in~\cite{bp}. This is because of two reasons. First, we do not weight the strata in the formula with the inverse order of the authomorphism group of their generic point. Second, we consider each possible enumeration of marked points only once, without multiplicities. We refer to~\cite{bp} for the explanation of the conventions that we do not keep here.

\subsection{Differential equation}

As we have already explained in the introduction, the Belorousski-Pandharipande relation gives
us a differential equation for $\Phi_0$, $\Phi_1$, $\Phi_2$, and $\Phi_2^{(1)}$. We illustrate this correspondence with examples.

If all variables are even, we have:
\begin{align}
\Delta_1 \rightsquigarrow\quad &
\frac{\d \Phi_{2}}{\d T_{0,i}}
\eta_{ij}
\frac{\d^3\Phi_0}{\d T_{0,j}\d T_{0,a}\d T_{0,k}}
\eta_{kl}
\frac{\d^3\Phi_0}{\d T_{0,l}\d T_{0,b}\d T_{0,c}} \\
& + 2\ terms\ obtained\ by\ permutations\ of\ \{a,b,c\}, \notag
\end{align}

\begin{align}
\Delta_2 \rightsquigarrow\quad &
\left(
\frac{\d \Phi_2^{(1)}}{\d T_{1,i}}
-
\frac{\d^2 \Phi_0}{\d T_{0,i} \d T_{0,k}}
\eta_{kl}
\frac{\d \Phi_2}{\d T_{0,l}}
\right)
\eta_{ij}
\frac{\d^4\Phi_0}{\d T_{0,j}\d T_{0,a}\d T_{0,b}\d T_{0,c}}  \label{m21} \\ \label{m22}
\Delta_3 \rightsquigarrow\quad &
\left(
\frac{\d^2 \Phi_2^{(1)}}{\d T_{1,a} \d T_{0,i}}
-
\frac{\d^2 \Phi_0}{\d T_{0,a} \d T_{0,k}}
\eta_{kl}
\frac{\d^2 \Phi_2}{\d T_{0,l} \d T_{0,i}}
\right)
\eta_{ij}
\frac{\d^3 \Phi_0}{\d T_{0,j}\d T_{0,b}\d T_{0,c}} \\
& + 2\ terms\ obtained\ by\ permutations\ of\ \{a,b,c\}, \notag
\end{align}
and so on. The metric $\eta_{ij}$ used here is given by the scalar product on $H_0$. We have:
$\eta_{ij}=(e_i,e_j)$, $\eta^{ij}=[\Pi_0]$.

Note that $\Delta_2$ is defined with the help of one $\psi$-class on $\oM_{2,1}$. Let $\pi\colon \oM_{2,n}\to \oM_{2,1}$ be the projection forgetting all marked points except for the first one. Then there is a formula relating $\psi_1$ and $\pi^*\psi_1$ on $\oM_{2,n}$. So, the differential expressions corresponding to these strata rely on this formula, which is exactly the first factor in Expression~\eqref{m21}. The same remark concerns $\Delta_2$, $\Delta_3$, and the pull-back of $\psi_1$ from $\oM_{2,2}$. In this case, the required formula is the first factor of Expression~\eqref{m22}.

\subsection{Theorem} We state our theorem.

\begin{theorem}\label{Belorousski-Pandharipande}
$\Phi_0$, $\Phi_1$, $\Phi_2$, and $\Phi_2^{(1)}$ satisfy the Belorousski-Pand\-ha\-ri\-pande relation.
\end{theorem}


\section{Proof}

In this section we prove Theorem~\ref{Belorousski-Pandharipande}. The proof is organized in two steps. First we consider the differential equation determined by the Belorousski-Pandharipande relation at zero point. It is proved by a straightforward calculation with tensors.

Then we can use the universal technique developed in~\cite{ls}. That is, for any differential equation proved at zero point by the same type of calculation as given below, we immediately obtain its proof at any point. This was done very carefully for the WDVV equation and less carefully for the Getzler relation in the last section of~\cite{ls}, and the argument for the Belorousski-Pandharipande relation is literally the same.

So, in subsections~\ref{relation-graphs}-\ref{other} we discuss only the simplest case of the Belorousski-Pandharipande relation. In subsection~\ref{relation-graphs} we rewrite it in terms of graphs; in subsection~\ref{outline} we explain what kind of calculation is to be done; in subsections~\ref{delta3} and~\ref{other} we give an example of such calculation for one stratum and discuss it for the other strata.

Finally, in subsection~\ref{reduction} we recall the basic idea from the last section of~\cite{ls} that completes the proof of any relation of this type.

\subsection{Relation in terms of graphs}\label{relation-graphs}

Consider the degree zero term of the power series obtained from $\Phi_0$, $\Phi_1$, $\Phi_2$, and $\Phi_2^{(1)}$ by the differential operator determined by $\Delta_i$. Slightly abusing notations, we denote it also by $\Delta_i$.
Then we have:
\begin{align}
\Delta_1 = & \frac{1}{16}\inspic{pictd.1}
+ \frac{1}{8}\inspic{pictd.2} + \frac{1}{8}\inspic{pictd.3} \\
\Delta_2 = & \frac{1}{12}\inspic{pictd.4}
+ \frac{1}{8}\inspic{pictd.5}
\end{align}
\begin{align}
\Delta_3=&\frac{1}{16}\inspic{pictd.6}+\frac{1}{8}\inspic{pictd.7}
+\frac{1}{8}\inspic{pictd.8}+\frac{1}{4}\inspic{pictd.9}\label{Del3}\\
&+\frac{1}{8}\inspic{pictd.10}+\frac{1}{8}\inspic{pictd.11}+\frac{1}{8}\inspic{pictd.12}\notag
\end{align}
and so on. We recall that a thick white point on an edge denotes the bivector $[\Pi_0]$ (or $[J\Pi_0]$).

\subsection{Outline of calculations}\label{outline}
We explain the proof of the simplest case of Theorem~1 . We have already expressed each $\Delta_i$ at zero point in terms of graphs with one or two edges marked by $[\Pi_0]$. Using the Leibniz rule for $Q$ and the $7$-term relation and $1/12$-axiom for $G_-$, we get out of $[\Pi_0]$ in our expressions. This way we obtain an expression of each $\Delta_i$ in terms of $60$ \emph{final graphs}. Then we substitute these expressions in Belorousski-Pandharipande relation~(\ref{BP-relation}), and we see that the coefficient of each final graph in this relation is equal to $0$. This proves the simplest case of our theorem.

Final graphs are listed in Appendix~A; final expressions of $\Delta_i$ are given in Appendix~B. Below, we explain how to get out of $[\Pi_0]$ in our graphs by the way of example (we give detailed calculations for $\Delta_3$).

\subsection{Calculations for $\Delta_3$}\label{delta3}

We consider the right hand side of Equation~(\ref{Del3}). Our goal is to get out of thick white points in these graphs. Finally, we must obtain an expression in terms of graphs from Appendix~A.

We carry out our calculations in two step. At the first step we consider each graph the right hand side of Equation~(\ref{Del3}) separately. At the second step we arrange the results of the first step and obtain an expression in final graphs.

\subsubsection{First step for the first picture}

We recall that $\Pi_0=\mathrm{Id}-QG_+-G_+Q$. Also we note that if we have an edge (not a loop) marked by $[\mathrm{Id}]$, then we can contract this edge. We have:
\begin{equation}
\frac{1}{16}\inspic{pictd.6}=\frac{1}{16}\inspic{calc.1}-
\frac{1}{16}\inspic{calc.2}-\frac{1}{16}\inspic{calc.3}
\end{equation}

We recall that $[Q,G_-G_+]=-G_-$ and $Qe_i=0$ for any $i$. Using these properties, the Leibniz rule for $Q$, and taking into account the symmetries of our graphs, we have:
\begin{align}
\frac{1}{16}\inspic{calc.2} & =\frac{1}{8}\inspic{calc.4}+\frac{1}{8}\inspic{calc.5} \\
\frac{1}{16}\inspic{calc.3} & =0
\end{align}

Therefore,
\begin{equation}\label{d3-1}
\frac{1}{16}\inspic{pictd.6}=\frac{1}{16}\inspic{calc.1}
-\frac{1}{8}\inspic{calc.4}-\frac{1}{8}\inspic{calc.5}
\end{equation}

\subsubsection{First step for all other pictures}

The same calculations for all other pictures give us:
\begin{align}
\frac{1}{8}\inspic{pictd.7}  = & \frac{1}{8}\inspic{calc.6} - \frac{1}{8}\inspic{calc.7}
- \frac{1}{8}\inspic{calc.8}\\
& -\frac{1}{8}\inspic{calc.9} -\frac{1}{8}\inspic{calc.10} \notag \\
\frac{1}{8}\inspic{pictd.8}  = & \frac{1}{8}\inspic{calc.11} - \frac{1}{4}\inspic{calc.12}
- \frac{1}{4}\inspic{calc.13} \\
\frac{1}{4}\inspic{pictd.9}  = & \frac{1}{4}\inspic{calc.14} - \frac{1}{4}\inspic{calc.15}
- \frac{1}{4}\inspic{calc.16} \\
& -\frac{1}{2} \inspic{calc.17} \notag
\end{align}

\begin{align}
\frac{1}{8}\inspic{pictd.10} = & \frac{1}{8}\inspic{calc.18} - \frac{1}{4}\inspic{calc.19}
- \frac{1}{8}\inspic{calc.20}
 - \frac{1}{8}\inspic{calc.21} \\
\frac{1}{8}\inspic{pictd.11} = & \frac{1}{8}\inspic{calc.22} - \frac{1}{4}\inspic{calc.23}
- \frac{1}{8}\inspic{calc.24} - \frac{1}{8}\inspic{calc.25} \\
\frac{1}{8}\inspic{pictd.12} = & \frac{1}{8}\inspic{calc.26} - \frac{1}{4}\inspic{calc.27}
- \frac{1}{8}\inspic{calc.28} - \frac{1}{8}\inspic{calc.29}\label{d3-7}
\end{align}

\subsubsection{Corollaries of $1/12$-axiom} In this section, we prove that some graphs in Equations~(\ref{d3-1})-(\ref{d3-7}) are equal to $0$.

\begin{lemma} Vector $\inspic{auux.1}$ is equal to $0$.
\end{lemma}

\begin{proof} Indeed, from $1/12$-axiom, it follows that
\begin{equation}
\inspic{auux.1}=\frac{1}{12}\inspic{auux.2}
\end{equation}
Since $G_-G_-G_+=0$, the last vector is equal to zero.
\end{proof}

From this lemma, it is obvious that
\begin{equation}
\inspic{calc.5}=\inspic{calc.8}=\inspic{calc.21}=\inspic{calc.25}=0.
\end{equation}

\begin{lemma} For any $i$ vector $\inspic{auux.3}$ is equal to $0$.
\end{lemma}

\begin{proof} First, we apply $1/12$-axiom, and then we apply the auxiliary lemma from~\cite{s}. We have:
\begin{equation}
\inspic{auux.3}=\frac{1}{12}\inspic{auux.4}=\frac{1}{12}\inspic{auux.5}+\frac{1}{12}\inspic{auux.6}
\end{equation}
Since $G_-G_-G_+=0$ and $G_-e_i=0$, the last two vectors are equal to zero.
\end{proof}

From this lemma, it is obvious that
\begin{equation}
\inspic{calc.7}=\inspic{calc.29}=0.
\end{equation}

\subsubsection{Corollaries of the $7$-term relation}

In this section, we list some corollaries of the $7$-term relation. We have:
\begin{multline}\label{17t}
\frac{1}{8}\inspic{calc.4}+\frac{1}{8}\inspic{calc.9}+\frac{1}{4}\inspic{calc.19}
+\frac{1}{4}\inspic{calc.23}\\
=\frac{1}{4}\inspic{cl.12}
\end{multline}

We prove this formula. For convenience, we split all these graphs in the same tensor pieces. We list the notations for these tensor pieces:
\begin{align}
x\otimes x & = \inspic{auux.7} & y &= \inspic{auux.8} & z &= \inspic{auux.9} & w &= \inspic{auux.10}
\end{align}
Note that $x,y,z$ are even vectors, but $w$ is an odd one. Equation~(\ref{17t}) is equivalent to
\begin{multline}\label{xyz-f}
\frac{1}{8}\int G_-(x^2) yzw+\frac{1}{8}\int G_-(x^2y)zw
+\frac{1}{4}\int G_-(xyz)xw
\\
+\frac{1}{4}\int G_-(xz)xyw=\frac{1}{4}\int G_-(x^2yz)w
\end{multline}
Also note that $G_-(x)=G_-(y)=G_-(z)=0$. Then from the $7$-term relation, it follows that
\begin{align}
G_-(x^2yz) & = 2G_-(xy)xz+2G_-(xz)xy+G_-(x^2)yz+G_-(yz)x^2 \label{rhs}\\
G_-(x^2y) & = 2G_-(xy)x+G_-(x^2)y \label{lhs1} \\
G_-(xyz) & = G_-(xy)z+G_-(xz)y+G_-(yz)x \label{lhs2}
\end{align}

Substituting Equations~(\ref{lhs1}) and~(\ref{lhs2}) in the left hand side of Equation~(\ref{xyz-f}), we get:
\begin{equation}
\int\left(
\frac{1}{4}G_-(x^2)yzw
+\frac{1}{2}G_-(xy)xzw
+\frac{1}{2}G_-(xz)xyw
+\frac{1}{4}G_-(yz)x^2w
\right)
\end{equation}
Substituting Equation~(\ref{rhs}) in the right hand side of Equation~(\ref{xyz-f}), we get exactly the same. This proves Equation~(\ref{xyz-f}) and therefore Equation~(\ref{17t}).

We prove in the same way that
\begin{align}
\frac{1}{8}\inspic{calc.10}+\frac{1}{4}\inspic{calc.27} &=\frac{1}{8}\inspic{cl.22} \\
\frac{1}{4}\inspic{calc.12}+\frac{1}{4}\inspic{calc.16} &=\frac{1}{6}\inspic{cl.13} \\
\frac{1}{4}\inspic{calc.13}+\frac{1}{8}\inspic{calc.24} &= 0
\end{align}

\begin{align}
\frac{1}{4}\inspic{calc.15} &= \frac{1}{12}\inspic{cl.23}
\end{align}

\begin{equation}\label{fDel3}
\frac{1}{2} \inspic{calc.17}+\frac{1}{8}\inspic{calc.20}+\frac{1}{8}\inspic{calc.28}=0
\end{equation}

\subsubsection{Final formula for $\Delta_3$}
Using Equations~(\ref{Del3})-(\ref{fDel3}) and notations for the final graphs from Appendix~A,
we get:
\begin{multline}
\Delta_3=\frac{1}{8}A_1+\frac{1}{16}A_2+\frac{1}{8}A_3+\frac{1}{8}F_1+\frac{1}{8}F_2
+\frac{1}{8}F_3+\frac{1}{4}F_4\\
-\frac{1}{4}D_2-\frac{1}{6}D_3-\frac{1}{8}H_1-\frac{1}{12}H_2
\end{multline}

\subsection{The other $\Delta_i$}\label{other} We carry out the same calculation for all other $\Delta_i$. If there are two thick white points in graphs for $\Delta_i$, then we get out of them successively. The results of these calculations are arranged in tables in Appendix~B.

If we substitute all these expressions for $\Delta_i$ in Belorousski-Pand\-ha\-ri\-pande relation, we identically get zero. This proves our theorem.

For a much more detailed exposition of our calculations, see~\cite{sn}.

\subsection{Reconstruction of the full proof}\label{reduction}

Now we explain what to do when parameters are not set to zero. In terms of graphs, this means
that for any $\Delta_i$ we are to consider graphs with the same structure as before but with an arbitrary number of additional leaves.

In~\cite{ls} the authors notice that these additional leaves can be gathered in some special operators. That is, instead of considering graphs with an arbitrary number of additional leaves, we can consider the same graphs as in the simplest case, but we substitute the vectors $E_0$ and $E_1$ on leaves and bivectors $[G_-G_+]$ and $[\Pi_0]$ on edges with new complicated vectors and bivectors.

These new vectors and bivectors depend on parameters and can be written down explicitely in terms of the Barannikov-Kontsevich solution of the Maurer-Cartan equation as it is done in~\cite{ls}.

Here is one subtlety related to the strata with one $\psi$-class. At this step we simultaneously have switched from
$\psi$-classes to a kind of pull-backs of $\psi$-classes. But is was proved in~\cite{s} that these pull-backs are related to $\psi$-classes via exactly the same formulas as in Gromov-Witten theory!

So, we take the same graphs as in the simplest case, we put new vectors on leaves and bivectors on edges, and we must prove \emph{exactly the same} relation as in the simplest case.

The main feature of this approach is that the properties of these new vectors and bivectors are almost the same as the properties of $E_0$, $E_1$, $[G_-G_+]$, and $[\Pi_0]$. So we can just repeat our argument with getting out of thick white points in graphs.

We refer to the last section of~\cite{ls} for the precise formulas for these new vectors and bivectors and lemmas describing their properties. In fact, this reconstruction of the full proof works for a rather large class of differential equations in cH-algebras, in particular all possible relations coming from geometry of the moduli space of curves are definitely in this class.


\appendix


\section{Final graphs}

\begin{align*}
A_1 &=\inspic{cl.1} & A_2 &=\inspic{cl.2} & A_3 &=\inspic{cl.3} & B_1 &=\inspic{cl.4} \\
B_2 &=\inspic{cl.5} & B_3 &=\inspic{cl.6} & C_1 &=\inspic{cl.7} & C_2 &=\inspic{cl.8} \\
C_3 &=\inspic{cl.9} & C_4 &=\inspic{cl.10} & D_1 &=\inspic{cl.11} & D_2 &=\inspic{cl.12} \\
D_3 &=\inspic{cl.13} & E_1 &=\inspic{cl.14} & E_2 &=\inspic{cl.15} & F_1 &=\inspic{cl.16} \\
\end{align*}

\begin{align*}
F_2 &=\inspic{cl.17} & F_3 &=\inspic{cl.18} & F_4 &=\inspic{cl.19} & G_1 &=\inspic{cl.20} \\
G_2 &=\inspic{cl.21} & H_1 &=\inspic{cl.22} & H_2 &=\inspic{cl.23} & H_3 &=\inspic{cl.24} \\
I_1 &=\inspic{cl.25} & I_2 &=\inspic{cl.26} & J_1 &=\inspic{cl.27} & J_2 &=\inspic{cl.28} \\
J_3 &=\inspic{cl.29} & K_1 &=\inspic{cl.30} & K_2 &=\inspic{cl.31} & L_1 &=\inspic{cl.32} \\
L_2 &=\inspic{cl.33} & M_1 &=\inspic{cl.34} & M_2 &=\inspic{cl.35} & M_3 &=\inspic{cl.36} \\
M_4 &=\inspic{cl.37} & M_5 &=\inspic{cl.38} & N_1 &=\inspic{cl.39} & N_2 &=\inspic{cl.40}
\end{align*}

\begin{align*}
N_3 &=\inspic{cl.41} & N_4 &=\inspic{cl.42} & N_5 &=\inspic{cl.43} & O_1 &=\inspic{cl.44} \\
O_2 &=\inspic{cl.45} & O_3 &=\inspic{cl.46} & O_4 &=\inspic{cl.47} & P_1 &=\inspic{cl.48} \\
P_2 &=\inspic{cl.49} & P_3 &=\inspic{cl.50} & Q_1 &=\inspic{cl.51} & Q_2 &=\inspic{cl.52} \\
Q_3 &=\inspic{cl.53} & R &=\inspic{cl.54} & S_1 &=\inspic{cl.55} & S_2 &=\inspic{cl.56} \\
S_3 &=\inspic{cl.57} & T_1 &=\inspic{cl.58} & T_2 &=\inspic{cl.59} & U &=\inspic{cl.60}
\end{align*}

\section{Results of calculations}

\begin{tabular}{|c|cccccccccc}

\hline
       & $A_1$ & $A_2$ & $A_3$ & $B_1$ & $B_2$ & $B_3$ & $C_1$ & $C_2$ & $C_3$ & $C_4$ \\
\hline
$\Delta_1$     &$\frac{1}{8}$ &$\frac{1}{16}$ &$\frac{1}{8}$ &$-\frac{1}{8}$ &$-\frac{1}{12}$ &$0$ &$0$ &$0$ &$0$ &$0$ \\
$\Delta_2$     &$0$ &$0$ &$0$ &$-\frac{1}{24}$ &$-\frac{1}{36}$ &$0$ &$0$ &$0$ &$0$ &$0$ \\
$\Delta_3$     &$\frac{1}{8}$ &$\frac{1}{16}$ &$\frac{1}{8}$ &$0$ &$0$ &$0$ &$0$ &$0$ &$0$ &$0$ \\
$\Delta_4$     &$\frac{1}{8}$ &$\frac{1}{16}$ &$\frac{1}{8}$ &$0$ &$0$ &$0$ &$\frac{1}{8}$ &$\frac{1}{8}$ &$\frac{1}{4}$ &$\frac{1}{8}$ \\
$\Delta_5$     &$0$ &$\frac{1}{24}$ &$0$ &$0$ &$0$ &$0$ &$0$ &$0$ &$0$ &$0$ \\
$\Delta_6$     &$0$ &$-\frac{1}{8}$ &$0$ &$0$ &$0$ &$0$ &$-\frac{1}{8}$ &$-\frac{1}{8}$ &$0$  &$-\frac{1}{8}$  \\
$\Delta_7$     &$0$ &$0$ &$0$ &$0$ &$0$ &$0$ &$\frac{1}{4}$ &$\frac{1}{4}$ &$0$ &$0$ \\
$\Delta_8$     &$0$ &$\frac{1}{8}$  &$0$ &$0$ &$0$ &$0$ &$0$ &$0$ &$0$ &$\frac{1}{4}$ \\
$\Delta_9$     &$0$ &$0$ &$0$ &$0$ &$0$ &$-\frac{1}{24}$ &$0$ &$0$ &$0$ &$0$ \\
$\Delta_{10}$  &$0$ &$\frac{1}{8}$ &$0$ &$0$ &$0$ &$0$ &$-\frac{1}{8}$ &$-\frac{1}{8}$ &$0$ &$\frac{1}{8}$ \\
$\Delta_{11}$  &$0$ &$\frac{1}{8}$ &$0$ &$0$ &$0$ &$\frac{1}{8}$ &$\frac{1}{8}$ &$\frac{1}{8}$ &$0$ &$-\frac{1}{8}$ \\
$\Delta_{12}$  &$0$ &$0$ &$0$ &$0$ &$0$ &$0$ &$0$ &$0$ &$0$ &$0$ \\
$\Delta_{13}$  &$0$ &$0$ &$0$ &$0$ &$0$ &$0$ &$0$ &$0$ &$0$ &$0$ \\
$\Delta_{14}$  &$0$ &$0$ &$0$ &$0$ &$0$ &$0$ &$0$ &$0$ &$0$ &$0$ \\
$\Delta_{15}$  &$0$ &$0$ &$0$ &$0$ &$0$ &$0$ &$0$ &$0$ &$0$ &$0$ \\
$\Delta_{16}$  &$\frac{1}{12}$ &$0$ &$\frac{1}{12}$ &$0$ &$0$ &$0$ &$0$ &$0$ &$0$ &$\frac{1}{4}$ \\
$\Delta_{17}$  &$-\frac{1}{4}$ &$0$ &$-\frac{1}{4}$ &$0$ &$0$ &$0$ &$0$ &$0$ &$-\frac{3}{4}$ &$-\frac{3}{4}$  \\
$\Delta_{18}$  &$\frac{1}{4}$ &$0$ &$\frac{1}{4}$ &$0$ &$0$ &$0$ &$0$ &$0$ &$\frac{3}{2}$ &$\frac{3}{4}$ \\
$\Delta_{19}$  &$0$ &$0$ &$0$ &$0$ &$0$ &$\frac{1}{12}$ &$0$ &$0$ &$0$ &$0$ \\
$\Delta_{20}$  &$0$ &$0$ &$0$ &$0$ &$0$ &$-\frac{1}{4}$ &$0$ &$0$ &$0$ &$0$ \\


 & & & & & & & & & & \\
\hline
               & $D_1$ & $D_2$ & $D_3$ & $E_1$ & $E_2$ & $F_1$ & $F_2$ & $F_3$ & $F_4$ & $G_1$ \\
\hline
$\Delta_1$     &$0$ &$0$ &$0$ &$0$ &$0$ &$0$ &$0$ &$0$ &$0$ &$0$ \\
$\Delta_2$     &$0$ &$\frac{1}{8}$ &$\frac{1}{12}$ &$0$ &$-\frac{1}{16}$ &$0$ &$0$ &$0$ &$0$ &$0$ \\
$\Delta_3$     &$0$ &$-\frac{1}{4}$ &$-\frac{1}{6}$ &$0$ &$0$ &$\frac{1}{8}$ &$\frac{1}{8}$ &$\frac{1}{8}$ &$\frac{1}{4}$ &$0$ \\
$\Delta_4$     &$0$ &$-\frac{1}{8}$ &$-\frac{1}{12}$ &$-\frac{1}{8}$ &$0$ &$0$ &$0$ &$0$ &$0$ &$0$ \\
$\Delta_5$     &$0$ &$-\frac{1}{8}$ &$0$ &$0$ &$\frac{1}{16}$ &$0$ &$0$ &$-\frac{1}{8}$ &$0$ &$0$ \\
$\Delta_6$     &$0$ &$\frac{1}{8}$ &$0$ &$\frac{1}{8}$ &$0$ &$0$ &$0$ &$\frac{3}{8}$ &$0$ &$0$ \\
$\Delta_7$     &$\frac{1}{8}$ &$\frac{1}{8}$ &$0$ &$0$ &$0$ &$0$ &$0$ &$-\frac{1}{8}$ &$0$ &$0$ \\
$\Delta_8$     &$0$ &$0$ &$0$ &$0$ &$0$ &$0$ &$0$ &$-\frac{1}{4}$ &$0$ &$0$ \\
$\Delta_9$     &$0$ &$0$ &$0$ &$0$ &$-\frac{1}{8}$ &$0$ &$0$ &$0$ &$0$ &$0$ \\
$\Delta_{10}$  &$0$ &$-\frac{1}{8}$ &$0$ &$-\frac{1}{8}$ &$0$ &$-\frac{1}{8}$ &$\frac{1}{8}$ &$-\frac{1}{8}$ &$0$ &$\frac{1}{8}$ \\
$\Delta_{11}$  &$0$ &$\frac{1}{8}$ &$0$ &$-\frac{1}{8}$ &$0$ &$\frac{1}{8}$ &$-\frac{1}{8}$ &$\frac{1}{8}$ &$0$ &$0$ \\
$\Delta_{12}$  &$0$ &$0$ &$0$ &$0$ &$0$ &$0$ &$0$ &$0$ &$0$ &$0$ \\
$\Delta_{13}$  &$0$ &$0$ &$0$ &$0$ &$0$ &$0$ &$0$ &$0$ &$0$ &$0$ \\
$\Delta_{14}$  &$0$ &$0$ &$0$ &$0$ &$0$ &$0$ &$0$ &$0$ &$0$ &$0$ \\
$\Delta_{15}$  &$0$ &$0$ &$0$ &$0$ &$0$ &$0$ &$0$ &$0$ &$0$ &$0$ \\
$\Delta_{16}$  &$-\frac{1}{4}$ &$0$ &$-\frac{1}{4}$ &$-\frac{1}{4}$ &$-\frac{1}{8}$ &$0$ &$-\frac{1}{4}$ &$0$ &$-\frac{1}{4}$ &$0$ \\
$\Delta_{17}$  &$\frac{1}{4}$ &$0$ &$\frac{1}{4}$ &$0$ &$0$ &$0$ &$\frac{3}{4}$ &$0$ &$\frac{3}{4}$ &$\frac{1}{4}$ \\
$\Delta_{18}$  &$0$ &$0$ &$\frac{1}{4}$  &$0$ &$0$ &$0$ &$-\frac{3}{4}$ &$0$ &$-\frac{3}{4}$   &$-\frac{1}{2}$ \\
$\Delta_{19}$  &$0$ &$0$ &$0$ &$0$ &$\frac{1}{4}$ &$0$ &$0$ &$0$ &$0$ &$0$ \\
$\Delta_{20}$  &$0$ &$0$ &$0$ &$\frac{1}{2}$ &$0$ &$0$ &$0$ &$0$ &$0$ &$-\frac{1}{4}$ \\
\end{tabular}

\begin{tabular}{|c|cccccccccc}
\hline
               & $G_2$ & $H_1$ & $H_2$ & $H_3$ & $I_1$ & $I_2$ & $J_1$ & $J_2$ & $J_3$ & $K_1$ \\
\hline
$\Delta_1$     &$0$ &$0$ &$0$ &$0$ &$0$ &$0$ &$0$ &$0$ &$0$ &$0$ \\
$\Delta_2$     &$0$ &$0$ &$0$ &$0$ &$0$ &$0$ &$0$ &$0$ &$0$ &$0$ \\
$\Delta_3$     &$0$ &$-\frac{1}{8}$ &$-\frac{1}{12}$ &$0$ &$0$ &$0$ &$0$ &$0$ &$0$ &$0$ \\
$\Delta_4$     &$0$ &$0$ &$0$ &$0$ &$0$ &$0$ &$0$ &$0$ &$0$ &$0$ \\
$\Delta_5$     &$0$ &$\frac{2}{16}$ &$0$ &$0$ &$0$ &$0$ &$0$ &$\frac{1}{8}$ &$0$ &$0$ \\
$\Delta_6$     &$0$ &$-\frac{1}{8}$ &$0$ &$0$ &$0$ &$0$ &$\frac{1}{4}$ &$-\frac{1}{4}$ &$0$ &$0$ \\
$\Delta_7$     &$0$ &$-\frac{1}{8}$ &$0$ &$0$ &$0$ &$\frac{1}{4}$ &$-\frac{1}{4}$ &$0$ &$0$ &$0$ \\
$\Delta_8$     &$0$ &$0$ &$0$ &$0$ &$\frac{1}{8}$ &$0$ &$-\frac{1}{4}$ &$\frac{1}{8}$ &$0$ &$0$ \\
$\Delta_9$     &$0$ &$0$ &$0$ &$\frac{1}{8}$ &$0$ &$0$ &$0$ &$0$ &$0$ &$\frac{1}{288}$ \\
$\Delta_{10}$  &$0$ &$\frac{1}{8}$ &$0$ &$0$ &$0$ &$0$ &$0$ &$0$ &$0$ &$0$ \\
$\Delta_{11}$  &$\frac{1}{4}$ &$-\frac{1}{8}$ &$0$ &$-\frac{1}{8}$ &$0$ &$0$ &$0$ &$0$ &$0$ &$-\frac{1}{96}$ \\
$\Delta_{12}$   &$0$ &$0$ &$0$ &$0$ &$0$ &$0$ &$0$ &$0$ &$0$ &$0$ \\
$\Delta_{13}$   &$0$ &$0$ &$0$ &$0$ &$0$ &$0$ &$0$ &$0$ &$0$ &$0$ \\
$\Delta_{14}$   &$0$ &$0$ &$0$ &$0$ &$0$ &$0$ &$0$ &$0$ &$0$ &$0$ \\
$\Delta_{15}$   &$0$ &$0$ &$0$ &$0$ &$0$ &$0$ &$0$ &$0$ &$0$ &$0$ \\
$\Delta_{16}$  &$\frac{1}{4}$ &$0$ &$\frac{1}{4}$ &$\frac{1}{4}$ &$0$ &$-\frac{1}{2}$ &$0$ &$0$ &$\frac{1}{4}$ &$0$ \\
$\Delta_{17}$  &$-\frac{1}{4}$ &$0$ &$-\frac{1}{4}$ &$-\frac{1}{4}$ &$-\frac{1}{2}$ &$\frac{1}{2}$ &$0$ &$0$ &$0$ &$0$ \\
$\Delta_{18}$  &$-\frac{1}{4}$ &$0$ &$-\frac{1}{4}$ &$-\frac{1}{4}$ &$\frac{3}{4}$ &$0$ &$0$ &$0$ &$-\frac{3}{4}$ &$0$ \\
$\Delta_{19}$  &$0$ &$0$ &$0$ &$-\frac{1}{4}$ &$0$ &$0$ &$0$ &$0$ &$0$ &$-\frac{1}{24}$ \\
$\Delta_{20}$  &$-\frac{1}{2}$ &$0$ &$0$ &$\frac{1}{4}$ &$0$ &$0$ &$0$ &$0$ &$0$ &$\frac{1}{8}$ \\
& & & & & & & & & & \\
%
%
\hline
               & $K_2$ & $L_1$ & $L_2$ & $M_1$ & $M_2$ & $M_3$ & $M_4$ & $M_5$ & $N_1$ & $N_2$ \\
\hline
$\Delta_1$     &$0$ &$0$ &$0$ &$0$ &$0$ &$0$ &$0$ &$0$ &$0$ &$0$ \\
$\Delta_2$     &$0$ &$0$ &$0$ &$0$ &$0$ &$-\frac{1}{96}$ &$0$ &$0$ &$0$ &$0$ \\
$\Delta_3$     &$0$ &$0$ &$0$ &$0$ &$0$ &$0$ &$0$ &$0$ &$0$ &$0$ \\
$\Delta_4$     &$0$ &$0$ &$0$ &$-\frac{1}{96}$ &$-\frac{1}{96}$ &$0$ &$0$ &$0$ &$0$ &$0$ \\
$\Delta_5$     &$0$ &$0$ &$0$ &$0$ &$0$ &$0$ &$0$ &$0$ &$0$ &$0$ \\
$\Delta_6$     &$0$ &$0$ &$0$ &$0$ &$0$ &$0$ &$0$ &$0$ &$\frac{1}{96}$ &$\frac{1}{96}$ \\
$\Delta_7$     &$0$ &$0$ &$0$ &$0$ &$0$ &$0$ &$0$ &$0$ &$-\frac{1}{48}$ &$-\frac{1}{48}$ \\
$\Delta_8$     &$0$ &$0$ &$0$ &$0$ &$0$ &$0$ &$0$ &$0$ &$0$ &$0$ \\
$\Delta_9$     &$\frac{1}{288}$ &$0$ &$0$ &$0$ &$0$ &$\frac{1}{96}$ &$0$ &$0$ &$0$ &$0$ \\
$\Delta_{10}$  &$0$ &$0$ &$0$ &$\frac{1}{96}$ &$\frac{1}{96}$ &$0$ &$0$ &$0$ &$0$ &$0$ \\
$\Delta_{11}$  &$-\frac{1}{96}$ &$0$ &$0$ &$0$ &$0$ &$0$ &$0$ &$0$ &$-\frac{1}{96}$ &$-\frac{1}{96}$ \\
$\Delta_{12}$  &$0$ &$0$ &$0$ &$0$ &$0$ &$-\frac{1}{8}$ &$0$ &$0$ &$0$ &$0$ \\
$\Delta_{13}$  &$0$ &$0$ &$0$ &$-\frac{1}{8}$ &$-\frac{1}{8}$ &$0$ &$0$ &$0$ &$-\frac{1}{8}$ &$-\frac{1}{8}$ \\
$\Delta_{14}$  &$0$ &$-\frac{1}{8}$ &$-\frac{1}{8}$ &$0$ &$0$ &$0$ &$0$ &$0$ &$\frac{1}{4}$ &$\frac{1}{4}$ \\
$\Delta_{15}$  &$0$ &$\frac{1}{8}$ &$\frac{1}{8}$ &$\frac{1}{8}$ &$\frac{1}{8}$ &$\frac{1}{8}$ &$\frac{1}{8}$ &$\frac{1}{4}$ &$-\frac{1}{8}$ &$-\frac{1}{8} $\\
$\Delta_{16}$  &$0$ &$0$ &$0$ &$0$ &$0$ &$\frac{1}{8}$ &$\frac{1}{96}$ &$\frac{1}{48}$ &$0$ &$0$ \\
$\Delta_{17}$  &$0$ &$0$ &$0$ &$\frac{1}{8}$ &$\frac{1}{8}$ &$0$ &$0$ &$0$ &$0$ &$0$ \\
$\Delta_{18}$  &$0$ &$\frac{1}{8}$ &$\frac{1}{8}$ &$0$ &$0$ &$0$ &$0$ &$0$ &$0$ &$0$ \\
$\Delta_{19}$  &$-\frac{1}{24}$ &$0$ &$0$ &$0$ &$0$ &$-\frac{1}{8}$ &$0$ &$0$ &$0$ &$0$ \\
$\Delta_{20}$  &$\frac{1}{8}$ &$0$ &$0$ &$-\frac{1}{8}$ &$-\frac{1}{8}$ &$0$ &$0$ &$0$ &$\frac{1}{8}$ &$\frac{1}{8}$ \\
\end{tabular}

\begin{tabular}{|c|cccccccccc}
\hline
               & $N_3$ & $N_4$ & $N_5$ & $O_1$ & $O_2$ & $O_3$ & $O_4$ & $P_1$ & $P_2$ & $P_3$ \\
\hline
$\Delta_1$      &$0$ &$0$ &$0$ &$0$ &$0$ &$0$ &$0$ &$0$ &$0$ &$0$ \\
$\Delta_2$      &$0$ &$0$ &$0$ &$0$ &$0$ &$0$ &$0$ &$0$ &$0$ &$0$ \\
$\Delta_3$      &$0$ &$0$ &$0$ &$0$ &$0$ &$0$ &$0$ &$0$ &$0$ &$0$ \\
$\Delta_4$      &$0$ &$0$ &$0$ &$0$ &$0$ &$0$ &$0$ &$0$ &$0$ &$0$ \\
$\Delta_5$     &$\frac{1}{96}$ &$0$ &$0$ &$-\frac{1}{96}$ &$\frac{1}{96}$ &$0$ &$0$ &$0$ &$-\frac{1}{96}$ &$-\frac{1}{96}$ \\
$\Delta_6$     &$-\frac{1}{48}$ &$0$ &$0$ &$\frac{1}{96}$ &$-\frac{1}{96}$ &$0$ &$0$ &$-\frac{1}{96}$ &$0$ &$0$ \\
$\Delta_7$     &$\frac{1}{96}$ &$0$ &$0$ &$0$ &$-\frac{1}{96}$ &$-\frac{1}{96}$ &$0$ &$0$ &$0$ &$\frac{1}{96}$ \\
$\Delta_8$      &$0$ &$0$ &$0$ &$0$ &$0$ &$0$ &$0$ &$0$ &$0$ &$0$ \\
$\Delta_9$     &$0$ &$0$ &$0$ &$0$ &$-\frac{1}{96}$ &$0$ &$-\frac{1}{96}$ &$0$ &$\frac{1}{96}$ &$0$ \\
$\Delta_{10}$   &$0$ &$0$ &$0$ &$0$ &$0$ &$0$ &$0$ &$0$ &$0$ &$0$ \\
$\Delta_{11}$  &$-\frac{1}{96}$ &$-\frac{1}{96}$ &$-\frac{1}{48}$ &$0$ &$\frac{1}{96}$ &$0$ &$\frac{1}{96}$ &$\frac{1}{96}$ &$0$ &$\frac{1}{96}$ \\
$\Delta_{12}$  &$-\frac{1}{8}$ &$0$ &$0$ &$\frac{1}{8}$ &$-\frac{1}{8}$ &$0$ &$0$ &$0$ &$\frac{1}{8}$ &$\frac{1}{8}$ \\
$\Delta_{13}$  &$\frac{1}{4}$ &$0$ &$0$ &$-\frac{1}{8}$ &$\frac{1}{8}$ &$0$ &$0$ &$\frac{1}{8}$ &$0$ &$0$ \\
$\Delta_{14}$  &$-\frac{1}{8}$ &$0$ &$0$ &$0$ &$\frac{1}{8}$ &$\frac{1}{8}$ &$0$ &$0$ &$0$ &$-\frac{1}{8}$ \\
$\Delta_{15}$  &$0$ &$0$ &$0$ &$0$ &$-\frac{1}{8}$ &$-\frac{1}{8}$ &$0$ &$-\frac{1}{8}$ &$-\frac{1}{8}$ &$0$ \\
$\Delta_{16}$   &$0$ &$0$ &$0$ &$0$ &$0$ &$0$ &$0$ &$0$ &$0$ &$0$ \\
$\Delta_{17}$   &$0$ &$0$ &$0$ &$0$ &$0$ &$0$ &$0$ &$0$ &$0$ &$0$ \\
$\Delta_{18}$   &$0$ &$0$ &$0$ &$0$ &$0$ &$0$ &$0$ &$0$ &$0$ &$0$ \\
$\Delta_{19}$  &0 &$0$ &$0$ &$0$ &$\frac{1}{8}$ &$0$ &$\frac{1}{8}$ &$0$ &$-\frac{1}{8}$ &$0$ \\
$\Delta_{20}$  &$\frac{1}{8}$ &$\frac{1}{8}$ &$\frac{1}{4}$ &$0$ &$-\frac{1}{8}$ &$0$ &$-\frac{1}{8}$ &$-\frac{1}{8}$ &$0$ &$-\frac{1}{8}$ \\
& & & & & & & & & & \\
%
%
\hline
               & $Q_1$ & $Q_2$ & $Q_3$ & $R$   & $S_1$ & $S_2$ & $S_3$ & $T_1$ & $T_2$ & $U$   \\
\hline
$\Delta_1$      &$0$ &$0$ &$0$ &$0$ &$0$ &$0$ &$0$ &$0$ &$0$ &$0$ \\
$\Delta_2$      &$0$ &$0$ &$0$ &$0$ &$0$ &$0$ &$0$ &$0$ &$0$ &$0$ \\
$\Delta_3$      &$0$ &$0$ &$0$ &$0$ &$0$ &$0$ &$0$ &$0$ &$0$ &$0$ \\
$\Delta_4$      &$0$ &$0$ &$0$ &$0$ &$0$ &$0$ &$0$ &$0$ &$0$ &$0$ \\
$\Delta_5$     &$0$ &$0$ &$0$ &$0$ &$-\frac{1}{48}$  &$0$ &$0$ &$\frac{1}{2304}$  &$\frac{1}{2304}$  &$\frac{1}{1152}$  \\
$\Delta_6$     &$0$ &$0$ &$0$ &$0$ &$\frac{1}{48}$  &$-\frac{1}{48}$  &$0$ &$0$ &$0$ &$0$ \\
$\Delta_7$     &$0$ &$0$ &$0$ &$0$ &$0$ &$\frac{1}{48}$  &$-\frac{1}{48}$  &$0$ &$0$ &$0$ \\
$\Delta_8$      &$0$ &$0$ &$0$ &$0$ &$0$ &$0$ &$0$ &$0$ &$0$ &$0$ \\
$\Delta_9$      &$0$ &$0$ &$0$ &$0$ &$0$ &$0$ &$0$ &$0$ &$0$ &$0$ \\
$\Delta_{10}$   &$0$ &$0$ &$0$ &$0$ &$0$ &$0$ &$0$ &$0$ &$0$ &$0$ \\
$\Delta_{11}$   &$0$ &$0$ &$0$ &$0$ &$0$ &$0$ &$0$ &$0$ &$0$ &$0$ \\
$\Delta_{12}$  &$0$ &$0$ &$0$ &$0$ &$\frac{1}{4}$  &$0$ &$0$ &$-\frac{1}{96}$  &$-\frac{1}{48}$ &$-\frac{1}{48}$  \\
$\Delta_{13}$  &$0$ &$0$ &$0$ &$0$ &$-\frac{1}{4}$  &$\frac{1}{4}$  &$0$ &$0$  &$0$ &$0$\\
$\Delta_{14}$  &$-\frac{1}{4}$ &$-\frac{1}{4}$ &$-\frac{1}{4}$ &$0$ &$0$ &$-\frac{1}{4}$  &$\frac{1}{4}$  &$0$ &$0$ &$0$ \\
$\Delta_{15}$  &$\frac{1}{4}$  &$\frac{1}{4}$  &$\frac{1}{4}$  &$\frac{1}{4}$  &$0$ &$0$ &$-\frac{1}{4}$  &$0$ &$0$ &$0$ \\
$\Delta_{16}$  &$0$ &$0$ &$0$ &$\frac{1}{48}$  &$0$ &$0$ &$0$ &$0$ &$0$ &$0$ \\
$\Delta_{17}$   &$0$ &$0$ &$0$ &$0$ &$0$ &$0$ &$0$ &$0$ &$0$ &$0$ \\
$\Delta_{18}$  &$\frac{1}{4}$  &$\frac{1}{4}$  &$\frac{1}{4}$  &$0$ &$0$ &$0$ &$0$ &$0$ &$0$ &$0$ \\
$\Delta_{19}$   &$0$ &$0$ &$0$ &$0$ &$0$ &$0$ &$0$ &$0$ &$0$ &$0$ \\
$\Delta_{20}$   &$0$ &$0$ &$0$ &$0$ &$0$ &$0$ &$0$ &$0$ &$0$ &$0$ \\
\end{tabular}

\end{document}